\documentclass[a4paper, 12pt]{article}

\usepackage[sort&compress]{natbib}
\bibpunct{(}{)}{;}{a}{}{,} 

\usepackage{amsthm, amsmath, amssymb, mathrsfs, multirow, url, subfigure}
\usepackage{graphicx} 
\usepackage{ifthen} 
\usepackage{amsfonts}
\usepackage[usenames]{color}
\usepackage{fullpage}
\usepackage{tikz}

\theoremstyle{plain}

\newtheorem*{thm0}{Theorem}
\newtheorem*{cor0}{Corollary}
\newtheorem*{fct}{False confidence theorem}

\theoremstyle{definition}

\theoremstyle{remark}

\newtheorem{ex}{Example}

\newcommand{\prob}{\mathsf{P}} 
\newcommand{\E}{\mathsf{E}}

\newcommand{\unif}{{\sf Unif}}
\newcommand{\nm}{{\sf N}}

\newcommand{\chisq}{{\sf ChiSq}}
\newcommand{\bet}{{\sf Beta}}

\newcommand{\RR}{\mathbb{R}}

\newcommand{\FF}{\mathbb{F}}
\newcommand{\MM}{\mathbb{M}}
\newcommand{\XX}{\mathbb{X}}

\newcommand{\TT}{\mathbb{T}}

\newcommand{\rlik}{\mathscr{P}}

\newcommand{\lPi}{\underline{\Pi}}
\newcommand{\uPi}{\overline{\Pi}}

\title{A possibility-theoretic solution to Basu's Bayesian--frequentist via media}
\author{Ryan Martin\footnote{Department of Statistics, North Carolina State University, {\tt rgmarti3@ncsu.edu}}}
\date{\today}

\begin{document}

\maketitle 

\begin{abstract}
Basu's {\em via media} is what he referred to as the middle road between the Bayesian and frequentist poles.  He seemed skeptical that a suitable {\em via media} could be found, but I disagree.  My basic claim is that the likelihood alone can't reliably support probabilistic inference, and I justify this by considering a technical trap that Basu stepped in concerning interpretation of the likelihood.  While reliable probabilistic inference is out of reach, it turns out that reliable {\em possibilistic} inference is not.  I lay out my proposed possibility-theoretic solution to Basu's {\em via media} and I investigate how the flexibility afforded by my imprecise-probabilistic solution can be leveraged to achieve the likelihood principle (or something close to it). 

\smallskip

\emph{Keywords and phrases:} conditional inference; fiducial argument; imprecise probability; inferential model; likelihood principle; validity.
\end{abstract}

\section{Introduction}
\label{S:intro}

Debabrata Basu (1924--2001) was a giant in the field who made fundamental contributions that have inspired generations of statisticians and helped shape the very core of our subject.  It's sincerely an honor and a privilege to contribute this manuscript for possible inclusion in the special volume of {\em Sankhya} in honor of Basu's birth centenary.

As I was thinking about what to contribute for this volume, I went back and reread some of Basu's classic works.  Of course, I read many of these papers when I was younger, even as a PhD student---my introductory stat theory class was taught by Anirban DasGupta who prompted us to {\em ``read Basu''} \citep[][p.~xvi]{basu.collection}---but I lacked the maturity to fully grasp their depth and quality at the time.  Now that I'm more experienced, I can better appreciate Basu's clarity and precision, along with his masterfully constructed examples.  Beyond that, I have the context to recognize the courage Basu had to critically challenge the leaders of both Fisher's and Neyman--Pearson--Wald's schools of thought.  Drawing inspiration from Basu's courage, here I make some similarly bold claims that I hope will stimulate discussions and help solidify our subject's foundations.

The title of this article references Basu's {\em via media}, a Latin phrase for ``the middle road.'' This comes from a remark he made in concluding his reply to the discussion of his monumental 1975 essay:
\begin{quote}
{\em The Bayesian and Neyman--Pearson--Wald theories of data analysis are the two poles in current statistical thought.  Today I find assembled before me a number of eminent statisticians who are looking for a {via media} between the two poles.  I can only wish you success in an endeavor in which the redoubtable R.~A.~Fisher failed.} \citep[][p.~269]{basu1975}
\end{quote}
The endeavor of Fisher's that Basu is referring to is, of course, the {\em fiducial argument};\footnote{One could argue that Basu is referring to Fisher's ideas on conditional inference. To me, however, Fisher's conditioning (like sufficiency) was largely motivated by a need to reduce a problem's dimension, without loss of information, so that his fiducial argument could be applied.} see, e.g., \citet{fisher1933, fisher1935a, fisher1973}, \citet{neyman1941}, \citet{lindley1958}, \citet{fraser1961b, fraser1965}, \citet{seidenfeld1992}, \citet{barnard1995}, \citet{dawid2020}, and many others, including the generalized fiducial developments summarized in, e.g., \citet{hannig2009} and \citet{hannig.review}.  Zabell and others have described fiducial inference as a sort of middle way:
\begin{quote}
{\em Fisher's attempt to steer a path between the Scylla of unconditional behaviorist methods which disavow any attempt at ``inference'' and the Charybdis of subjectivism in science was founded on important concerns, and his personal failure to arrive at a satisfactory solution to the problem means only that the problem remains unsolved, not that it does not exist}. \citep[][p.~382]{zabell1992}
\end{quote}
In addition to the obvious similarities in how the two authors characterize Fisher's efforts to strike this balance, they both leave open the {\em possibility} (pun intended) of a resolution, though Basu's remark falls short of Zabell's on the optimism scale.  To me, it's imperative to the long-term success of the field of statistics that the {\em via media} be found and, fortunately, a solution is currently available.  The high-level goal of this paper is to motivate and explain this solution, while making connections to Basu's work.

Part of motivating this solution---and even why a solution is needed---is understanding how the priorities of today's statisticians differ from those in Basu's time, when our subject was taking shape.  
Right or wrong, many view statisticians' role in science as connecting the scientific problem~P to a suitable statistical method~M to apply.  In fact, the courses offered to non-statistics students tend to focus on enumerating the standard $(\text{P}, \text{M})$ pairs, e.g., comparing treatments via randomized experiments and the analysis of variance.  In the old days, when a scientist encountered a new or unfamiliar problem $\text{P}'$, he'd probably consult with a statistician.  Nowadays, the scientist can immediately find a method~$\text{M}'$ to apply in his problem $\text{P}'$ by consulting Google, no statistician required.  So, while there are exceptions, modern-day statisticians' involvement in the scientific process is more indirect, by having an article that appears near the top of Google's search results.  Consequently, those entertaining and expertly-crafted hypothetical dialogues between scientist and statistician in, e.g., \citet{basu1975, basu1980} and \citet{bergerwolpert1984}, designed to shine light on our foundational questions/answers, no longer ring true.  

As I explain in Section~\ref{S:priorities} below, this shift in the role that statisticians play in the scientific process, from direct to indirect, marks a change in statisticians' priorities.  With respect to the two poles that Basu mentioned, now almost everyone is gathered around the frequentist pole---even the Bayesians!  The methods-developing statistician simply can't ignore frequentist considerations, the very same frequentist considerations that the aforementioned dialogues crafted by Basu and others aimed to show were irrelevant or downright silly.  This is an unacceptably wide gap between our practical priorities and what our current foundations say.  On the one hand, to jump to Basu's preferred Bayesian pole is tantamount to ignoring the modern priorities.  On the other hand, to stay gathered around the frequentist pole and ignore the foundational issues raised by Basu and others is tantamount to concluding that those insights were wrong and/or no longer relevant.  In both cases we end up losing our seat at the data science table: in the former, we're ignoring modern priorities and, in the later, we're admitting that our history and experience gives us no upper-hand over our new competitors.  Neither of these are desirable outcomes, so a {\em via media} is imperative for our field's long-term success.  

Section~\ref{S:evolve} presents how I expect the {\em via media} to look.  Like the Bayesian pole, it offers fixed-data ``probabilities'' which can be used for ``probabilistic reasoning'' and inference; like the frequentist pole, these ``probabilities'' satisfy a certain {\em validity} property which implies that the procedures derived from them have error rate guarantees.  This may sound too good to be true, and it is.  The catch is that what I referred to above as ``probabilities'' aren't probabilities in the familiar sense; they're {\em imprecise probabilities} or, more specifically, they're {\em possibilities}.  The shift from probability theory/calculus to the corresponding possibility theory/calculus is technically simple, but a fundamental change like this can be a conceptually large pill to swallow.  My claim is that the likelihood alone can't reliably support probabilistic inference, so sticking with probability isn't an option.  To justify this claim, in Section~\ref{S:likelihood}, I highlight a technical trap that Basu stepped in, related to the well-known fact that a likelihood function isn't a probability density, i.e., it has no inherent differential element.  Other authors \citep[e.g.,][]{shafer1982, wasserman1990b} have suggested that the likelihood is more appropriately processed as a differential element-free possibility contour, but their proposals don't go far enough. 

Having explained the high-level vision behind my proposed {\em via media} and justified the transition from probabilistic to possibilistic thinking, I provide a more detailed description of its implementation in Section~\ref{S:new}.  As I explain, the proposal shares some similarities with what's commonly done in statistical practice, but it's part of a framework that itself is very different.  This, to me, is exactly what we'd expect from a suitable {\em via media}---it must be different from the two poles, but not unrecognizably different.  I see this proposal as a modern, likelihood-centric version of the {\em inferential model} (IM) framework put forward in \citet{imbasics, imbook}; for many more details, see \citet{imchar, martin.partial, martin.partial2}.  A few Basu-spirited illustrations of my proposal are also presented here.  

In Section~\ref{S:favorites}, I consider one of Basu's favorites---the likelihood principle \citep{birnbaum1962}---and how the new perspectives afforded by the imprecise-probabilistic formulation of my proposed {\em via media} can be beneficial.  In particular, note that my basic proposal in Section~\ref{S:new} doesn't satisfy the likelihood principle, but it's clear how it can be made to do so without sacrificing on the solution's validity, provided that the data analyst is willing to give up some of their solution's efficiency.  It's also possible to be valid and {\em partially} satisfy the likelihood principle, e.g., to be valid with respect to some user-identified set of plausible stopping rules but not to others, thereby balancing both the efficiency and the stopping-rule invariance that are relevant around the frequentist and Bayesian poles, respectively.  The paper concludes with a brief discussion in Section~\ref{S:discuss}.

\section{Priorities have changed}
\label{S:priorities}

Despite the very powerful foundational arguments put forward by \citet{savage1972}, \citet{basu1975}, \citet{bergerwolpert1984}, and many others in support of a fully-conditional, likelihood-centric approach to statistical inference, it's fair to say that there's effectively no sign of this way of thinking in modern statistics research---even among Bayesians.  My claim is that statisticians' priorities have changed. 

There are exceptions, of course, but today's academic statisticians, for the reasons I explained above, are almost exclusively focused on the development of {\em statistical methods}, i.e., specific tools and software intended to be used off-the-shelf by scientists working on the scientific front lines.  The scientist is motivated by results, so their top priority is that a statistical method ``works.''  That is, they're not going to apply a method off-the-shelf unless it's been demonstrated to ``work'' in some meaningful sense.  This begs the question: in what sense could a method ``work'' that would be meaningful to scientists?  It seems necessary that the method has been demonstrated to give a ``right answer'' in most of the cases in which it's applied.  Then a scientist who believes that his problem is similar to those in which it's been demonstrated that the method typically gives a ``right answer'' has no reason to doubt that his application is one of those typical cases and, consequently, no reason to doubt the result of that method applied to his problem. 


The reader surely recognizes that my description of what it means for a method to ``work'' is very much frequentist.  The reader surely is also aware that these frequentist considerations, and my definition of ``works,'' don't align with the fully-conditional likelihood/Bayesian considerations of Basu and others.  ``Don't align'' is an understatement, these two considerations are almost completely incompatible---if a method ``works'' in the sense above, then it almost always falls short of Basu's foundational bar.  So where does this leave us?  For sure, a subject that's central to the advancement of knowledge shouldn't abandon its foundations altogether for the {\em priorit\'e du jour}.  But it's similarly embarrassing for the same subject to hold up a foundational standard that's not generally taken seriously by today's methods-developing statisticians.  Even modern Bayesian methods fail to meet Basu's high standard.  It's of course well-known that the (common) use of default priors is incompatible with the likelihood principle.  Moreover, there's a relevant selection bias in the Bayesian literature: the Bayesian methods that appear in publications tend to be those that have been demonstrated to ``work'' in the sense above, either theoretically or empirically.\footnote{An anonymous reviewer asked for some evidence to support this claim, so I scanned the 2022 volume of {\em Bayesian Analysis} and noted which papers offered a demonstration that the proposed method ``works'' based on a proof of consistency or a simulation study.  Of the 45 published papers, I counted 37 papers (about 82\%) with a significant focus on the proposed method's frequentist performance.} These demonstrations fix the data-generating process, so adopting a Bayesian solution in an application because the method ``works'' under certain data-generating processes is a violation of the likelihood principle.

%
%
%
%
%
%
%
%

To be fair, the Scylla at the frequentist pole isn't any more pleasant than the Charybdis at Bayesian pole.  Aside from not really addressing the question of {\em inference}, the pure performance focus hasn't proved to ensure reliability, as the replication crisis has revealed.  The take-away, again, is that the problem of statistical inference can't be fully resolved at either of the extreme poles; the nuance if a genuine {\em via media} is necessary.

\section{Towards a {\em via media}}
\label{S:evolve}

As Zabell wrote in the quote above, the fact that Fisher's attempts to find this {\em via media} failed doesn't mean that it can't be found.  But it will require some outside-the-box ideas, and I'll share these ideas in the subsequent sections.  First, I should explain how I think this {\em via media} ought to look, since what I have in mind is quite different from what is currently done by both Bayesians and frequentists.  My main thesis is as follows: 
\begin{itemize}
\item a fully satisfactory theory of statistical inference ought to produce reliable, data-dependent ``probabilities'' based on which probabilistic reasoning can be made, i.e., if the data-dependent ``probability'' assigned to a hypothesis is small, then that should provide good reason to doubt the truthfulness of the hypothesis;
\vspace{-2mm}
\item but the data alone can't reliably support the construction of data-dependent ``probabilities'' that are genuine probabilities like in Kolmogorov's theory;
\vspace{-2mm}
\item so, to achieve both the probabilistic reasoning that's advantageous for single-case, data-driven inference and the reliability that's necessary from today's perspective with a focus on methods-development, this {\em via media} can't be contained in the existing/standard theory of probability, i.e., the ``probabilities'' I'm referring to above can't literally be probabilities in the sense of Kolmogorov.  
\end{itemize} 
My justification for the claim in the second bullet point will come in the next section, where things get more technical.  In the remainder of this section, I want to focus on the first and third bullet points, which are more conceptual in nature.  

For the first bullet: why is probabilistic reasoning so important?  A common criticism of the frequentist theory of inference, which isn't based on probabilistic reasoning, is that significance levels, coverage probabilities, etc.~are pre-data calculations---they don't speak to what the observed data actually say about the unknown being inferred.  P-values aim to bring the observed data into the uncertainty quantification picture, but these too are often (unjustifiably) criticized because they're not probabilities, not measures of the strength of evidence, etc. 
More recently, some authors, especially Deborah Mayo, have been calling for more than what the classical frequentist solutions offer.  She argues in \citet{mayo.book.2018} that, in addition to determining if a hypothesis is incompatible with the data, via tests and confidence sets, it's important that scientists can ``probe'' deeper into those hypotheses that are compatible with the data to find sub-hypotheses that might actually be supported by data.  This probativeness feature comes fairly naturally when inference is based on data-dependent ``probabilities,'' but not otherwise; Mayo suggests supplementing the frequentist methods with a so-called {\em severity measure} designed to offer probativeness.  Suffice it to say that there are real, practical advantages to probabilistic reasoning that the classical frequentist solutions fail to offer, but these advantages don't come for free just by choosing to write down (artificial) probabilities.  

It's in the third bullet where the {\em via media} starts to reveal itself.  Recall that Basu's poles correspond to probability (Bayesian) and not-probability (frequentist).  From this perspective, it seems almost obvious that the middle-ground must somehow be both probability and not-probability simultaneously.  ``Fisher's biggest blunder'' \citep{efron1998} was just his failure to see that the {\em via media} can't be achieved entirely within the theory of probability.  What I/we have now that Fisher didn't have is (the benefit of hindsight and) more than 60 years of developments---starting with Art Dempster's seminal work in the 1960s \citep[e.g.,][]{dempster1966, dempster1967}---in the theory of {\em imprecise probability}.  What I'm proposing in Section~\ref{S:new} falls under the umbrella of imprecise probability, but it's both drastically different from Dempster's approach and surprisingly similar to ideas that can be found in standard statistics textbooks.  But before I can get into these details, I need to justify the claim in my second bullet point above, which I'll do next.

\section{Likelihood and inference}
\label{S:likelihood}

\subsection{What likelihood can't do}

To set some notation, let $X \in \XX$ denote the observable data and write $\prob_\Theta$ for the posited statistical model, depending on an unknown parameter $\Theta \in \TT$.  As is customary, here I'll write $\Theta$ for the unknown true parameter value, saving $\theta$ and $\vartheta$ to denote generic parameter values.  For fixed $\theta \in \TT$, suppose that $\prob_\theta$ admits a probability mass or density function $p_\theta$ on $\XX$, and define the likelihood function at the observed $X=x$, as $L_x(\theta) = p_\theta(x)$. 
The name ``likelihood'' was coined by Fisher and part of the motivation behind this choice of name was to emphasize that, notwithstanding the obvious connection between the likelihood function and the model's probability density/mass function, the likelihood is indeed fundamentally different from probability.  In particular:
\begin{quote}
{\em The function of the $\theta$'s... is not however a probability and does not obey the laws of probability; it involves no differential element $d\theta_1 \, d\theta_2 \, d\theta_3$...; it does none the less afford a rational basis for preferring some values of $\theta$, or combination of values of the $\theta$'s, to others.} \citep[][p.~552]{fisher1930}
\end{quote}
Despite the warnings, many have not taken this seriously---including Fisher himself!  Indeed, as \citet[][p.~33]{basu1975} points out, there are cases in which Fisher's fiducial argument, his proposed {\em via media}, produces a solution that's equivalent to treating the likelihood as if it were a probability density/mass function for $\Theta$, given $X=x$.  If the likelihood doesn't determine a probability distribution for $\Theta$, and if the fiducial argument can produce a solution that's a probability determined by the likelihood, then isn't that an obvious sign something's wrong with the fiducial argument itself? 

But Basu stepped into this trap too.  In \citet[][Sec.~8]{basu1975}, he proposes the construction of a data-dependent probability distribution for $\Theta$ on $\TT$ based on a normalized likelihood function,\footnote{Basu actually assumes $\TT$ is finite and defined the above expression with the integrals replaced by sums; see \eqref{eq:shackle}.  I'm writing integrals here only because it's more common for the parameter space to be a continuum; none of what I have to say here depends on this choice.} which, in the present notation, is 
\begin{equation}
\label{eq:lik.normal}
\bar L_x(A) = \frac{\int_A L_x(\theta) \, d\theta}{\int_\TT L_x(\theta) \, d\theta}, \quad A \subseteq \Theta. 
\end{equation}
Constructing a probability by suitably normalizing the likelihood function as above seems natural and, following a detailed analysis, \citet[][p.~33]{basu1975} concludes with: 
\begin{quote}
{\em The author can find no logical justification for the often repeated assertion that likelihood is only a point function and not a measure.  He does not see what inconsistencies can arise from [treating it as a measure].}
\end{quote}
Problems arise because, as Fisher emphasized, the likelihood has no differential element ``$d\theta$.''  While introducing ``$d\theta$'' and normalization via integration might seem innocuous, this isn't free of consequences.  Of course, $A \mapsto \bar L_x(A)$ is a measure, so it's additive and monotone: in particular, $A \subseteq B$ implies $\bar L_x(A) \leq \bar L_x(B)$.  That $A$ can't be more compatible with data than $B$ is perfectly logical, but it'll virtually always be that $\bar L_x(B)$ is {\em strictly greater} than $\bar L_x(A)$.  For example, suppose that $X \sim \nm(\Theta,1)$, and that $x=7$ is observed.  Consider two hypotheses about the unknown $\Theta$: $A=[7.7, 8]$ and $B=[7.7, 20]$.  Clearly, $A$ is a proper subinterval of $B$, and $B$ is much wider than $A$, which implies that $\bar L_x(B) \gg \bar L_x(A)$; in this particular case, $\bar L_x(B) \approx  3\bar L_x(A)$.  But is there any sense in which $B$ is {\em strictly more} compatible with the data than $A$?  No---it's obvious that the inclusion of points that are relatively incompatible with the data doesn't make the hypothesis more compatible with the data.  That's the point I think Fisher was trying to make when he emphasized the likelihood involves no differential element.  

Similar points were made by economist G.~L.~S.~Shackle in the mid-1900s.  Like in \citet[][p.~29]{basu1975}, Shackle had in mind a finite space $\TT$ and was entertaining the option of assigning plausibility\footnote{Shackle didn't specifically mention likelihood in his analysis, but my choice to make this point using likelihood is consistent with Shackle's remarks.} to individual elements as 
\begin{equation}
\label{eq:shackle}
\bar L_x(\theta) = \frac{L_x(\theta)}{\sum_{\vartheta \in \TT} L_x(\vartheta)}, \quad \theta \in \TT. 
\end{equation}
Note that the above relationship forces the mass assigned to an individual $\theta$ to depend on the cardinality of $\TT$.  Shackle argues emphatically that the size of (the hypothesis space) $\TT$ ought not to influence the plausibility of an individual (hypothesis) $\theta$.
\begin{quote}
{\em To allow the size of the crowd of hypotheses...
to influence the value of the [plausibility] assigned to any particular hypothesis, would be like weakening one's praise for the chief actors in a play on the ground that a large number of supers were also allowed to cross the stage.} \citep[][p.~51]{shackle1961}
\end{quote}

This begs a fundamental question: does introducing an artificial differential element detail have any practical consequences?  Yes!  It implies existence of true hypotheses $A$ for which $\bar L_X(A)$ tends to be small as a function of $X$ and, similarly, the existence of false hypotheses $B$ for which $\bar L_X(B)$ tends to be large as a function of $X$.  Since one would be inclined, e.g., to doubt the truthfulness of a hypothesis for which $\bar L_X(A)$ is small, this counter-intuitive behavior raises serious practical concerns about the reliability of inferences based on the normalized likelihood---which is very much relevant to the methods-developing statistician and to the scientist who uses these methods.  The root cause of this undesirable behavior is obvious: there are small sets that contain $\Theta$ and large sets that don't.  That is, the size itself of a hypothesis has no bearing on whether it's true or false and, therefore, no bearing on how compatible it is with the data.  This intuition is captured by the likelihood and its lack of a differential element.  But the integral-driven normalization forces additivity, which allows the (irrelevant) size of the hypothesis to become relevant, hence the undesirable behavior.  This is captured in the main result of \citet{balch.martin.ferson.2017}, though the connection here to Fisher's warning about the likelihood's lack of a differential element is apparently new. 

\begin{fct}[\citealt{balch.martin.ferson.2017}]
Let $\bar L_x(\cdot)$ be Basu's $x$-dependent probability distribution \eqref{eq:lik.normal} based on normalizing the likelihood via integration.  Then for any pair $(\lambda,\alpha) \in [0,1]^2$ and for any $\Theta \in \TT$, there exists a hypothesis $A \subset \TT$ such that 
\[ A \not\ni \Theta \quad \text{and} \quad \prob_\Theta\{ \bar L_X(A) > \lambda \} > \alpha. \]
The same conclusion holds true for any data-dependent probability measure, not just that in \eqref{eq:lik.normal}, including any Bayes posterior or fiducial distribution. 
\end{fct}

In words, the false confidence theorem states that there are false hypotheses to which the artificially additive posterior $\bar L_X(\cdot)$ tends to assign high probability; this creates a risk of systematically misleading conclusions and raises doubts about the reliability of inference.  This is my justification for the claim made in the second bullet point in Section~\ref{S:evolve} above: the likelihood (data + model) alone can't reliably support genuine data-dependent {\em probabilities} and the associated probabilistic inference.  

To be clear, the above issues are unrelated to the choice of dominating measure: one can't sidestep the difficulties raised by the false confidence theorem by introducing a default prior density/mass function before normalization.  The point, again, is that a hypothesis's size has no direct bearing on its compatibility with the data, and yet it's relevant to any Lebesgue integral.  A hypothesis's size is relevant only if the data analyst {\em believes} that it is, i.e., if he's willing to introduce a genuine prior distribution that connects size to credence.  That is, if the prior probabilities are real and part of the posited model, then the differential element is meaningful and there's no false confidence.  If the prior probabilities are artificial, then there are no guarantees: {\em [Bayes's formula] does not create real probabilities from hypothetical probabilities} \citep[][p.~249]{fraser.copss}.

\subsection{What likelihood can do}

If the likelihood doesn't have a differential element and, therefore, doesn't reliably support probabilistic inference, but there's still something that can be done!  The point is: probability theory is not the only uncertainty quantification game in town.  Starting with Dempster's seminal work in the 1960s, there have been major developments in what's called {\em imprecise probabilities}; see, e.g., the books by \citet{shafer1976}, \citet{dubois.prade.book}, \citet{walley1991}, \citet{lower.previsions.book}, \citet{imprecise.prob.book}, and \citet{cuzzolin.book}.  The simplest among the imprecise probability models is {\em possibility theory}, with close connections to fuzzy set theory \citep[e.g.,][]{zadeh1965, hanss.book}, which dates back to \citet{shackle1961} and \citet{zadeh1978}; much of modern possibility theory is based on \citet{dubois.prade.book} and the extensive subsequent work by the same authors.  It is a general-purpose theory of uncertainty quantification applied throughout science and engineering; statistical applications are discussed in \citet{dubois2006} and the references therein, but the perspective I share here and in the next section is new. 

A simple idea, similar to Basu's, starts by defining the {\em relative likelihood} function
\begin{equation}
\label{eq:rho}
\eta_x(\theta) = \frac{L_x(\theta)}{\sup_\vartheta L_x(\vartheta)}, \quad \theta \in \TT.
\end{equation}
Note the difference in normalization---supremum versus Basu's integration---so $\eta_x$ isn't a probability.  
But \eqref{eq:rho} is the driver behind the proposal in \citet{shafer1976, shafer1982}, which is developed further in \citet{wasserman1990b}, \cite{denoeux2014}, and elsewhere.  This determines a very special imprecise probability structure which has a few different names: here I adopt the possibility theory terminology, so I'll refer to $\theta \mapsto \eta_x(\theta)$ in \eqref{eq:rho} as a {\em possibility contour}.  Mathematically, what distinguishes $\eta_x$ as a possibility contour is that, first, it takes values in $[0,1]$ and, second, that it satisfies $\sup_\theta \eta_x(\theta) = 1$.  The contour's extension to a {\em possibility measure} supported on general hypotheses is defined as 
\[ \eta_x(A) = \sup_{\theta \in A} \eta_x(\theta), \quad A \subseteq \TT. \]
This is analogous to Basu's $\bar L_x$, just with possibility calculus\footnote{Note that possibility calculus can be described via a suitable Choquet integral instead of the familiar Lebesgue integral; see \citet{choquet1953} and \citet[][App.~C]{lower.previsions.book}.} instead of probability calculus.  This way of processing the likelihood function has a number of desirable properties, e.g., it's completely driven by the likelihood-based ranking of the parameter values so it doesn't require introduction of an artificial differential element.  

Of course, the set-function $A \mapsto \eta_x(A)$ isn't a measure in the usual sense, but it does have some similar properties.  In addition to $\eta_x(\cdot) \geq 0$ and $\eta_x(\TT) = 1$, the possibility measure is {\em maxitive}\footnote{Maxitive means $\eta_x( \bigcup_{k=1}^\infty A_k ) = \sup_{k \geq 1} \eta_x(A_k)$ for all $A_1,A_2,\ldots \subseteq \TT$.} 
which implies sub-additivity, in particular, 
\begin{equation}
\label{eq:non-additive}
1 \leq \eta_x(A) + \eta_x(A^c), \quad \text{for all $A \subseteq \TT$}. 
\end{equation}
Maxitivity also implies monotonicity, but not the kind of strict monotonicity that often holds for probabilities.  Reconsider the simple $X \sim \nm(\Theta,1)$ illustration above, with $A=[7.7, 8]$ and $B=[7.7, 20] = A \cup (8,20]$.  While Basu's $\bar L_x$ has $\bar L_x(A) \ll \bar L_x(B)$, the possibility measure has $\eta_x(A) = \eta_x(B)$, as one would expect:  the inclusion of an interval $(8,20]$ that contains ``less-likely'' values shouldn't increase the compatibility with $x$. 

Mathematics aside, since $\eta_x$ isn't a probability, we don't have access to the full power of probabilistic reasoning.  But a one-sided version is available, what I'll call here {\em possibilistic reasoning}.  That is, possibility theory allows for a direct refutation of a hypothesis ``$\Theta \in A$'' by showing that $\eta_x(A)$ is small.  However, unlike with probabilistic reasoning, if $\eta_x(A)$ is large, then that's not enough to conclude that there's support for ``$\Theta \in A$,'' since \eqref{eq:non-additive} doesn't rule out the case that both $\eta_x(A)$ and $\eta_x(A^c)$ are large.  In possibilistic reasoning, we need {\em both} $\eta_x(A)$ large and $\eta_x(A^c)$ small to find support for ``$\Theta \in A$.''  

What does it mean for $\eta_x(A)$ to be ``small'' or ``large''?  The methods-developing statistician must suggest to potential users how to make these judgments and, if his method is going to demonstrably ``work,'' then he similarly must take this small/large judgment seriously.  One can tailor these small/large possibility thresholds to the problem at hand, or perhaps rely on asymptotic theory to get some unification, but that's different from probabilistic reasoning.  Indeed, recall that probability has the same scale across every example to which it's applied: a numerical probability of 0.1 means the same thing whether it's the probability of rain tomorrow or the probability of a patient responding favorably to a new cancer treatment.  The basic likelihood-to-possibility setup presented above doesn't share this invariance, i.e., the small/large scale that ``works'' depends crucially on features of the application at hand.  But a different possibility-theoretic framework can do it, as I explain next.


\section{A possibility-theoretic {\em via media}}
\label{S:new}

The details below are simultaneously both new and familiar, i.e., there are close connections with classical theory but the possibility-theoretic details that make it a full-blown framework are recent developments and are likely unfamiliar to most readers.  For the sake of space, I'll only present the immediately-relevant aspects of this theory.  In particular, I'll not present the (arguably most interesting) details that showcase how the framework easily incorporates {\em partial prior information} about $\Theta$.  The partial-prior angle is crucial for at least two reasons: first, it's what creates new opportunities for improved methods and, second, it's what justifies this proposal as a bona fide {\em via media} between the Bayesian and non-Bayesian poles.  The interested reader can consult \citet{martin.partial, martin.partial2}. 

It turns out that the relative likelihood function is still very relevant here \citep[see][Sec.~5.1]{martin.partial2}.  But since its role is a bit different, I'm going to use a slightly modified notation: I'll write $\eta(x,\theta)$ instead of $\eta_x(\theta)$.  The key idea is that the likelihood offers a data-dependent partial order on $\TT$, but even the relative likelihood is lacking a universal scale that ``works'' for all applications.  Following the key developments in \citet{hose2022thesis}, my proposal is to calibrate the relative likelihood in a principled way to construct a new possibility contour---and corresponding possibility measure---that has the same partial order on $\TT$ as the likelihood but is universally scaled and provably ``works.''  Define this new likelihood-based possibility contour for $\Theta$ as 
\begin{equation}
\label{eq:pl}
\pi_x(\theta) = \prob_\theta\{ \eta(X,\theta) \leq \eta(x,\theta)\}, \quad \theta \in \TT. 
\end{equation}
The reader may recognize this as a sort of p-value determined by the relative likelihood.  This is for the special case where prior information about $\Theta$ is vacuous; if (partial) prior information is available, then a different possibility contour emerges.  In \citet{martin.partial2}, it's shown that \eqref{eq:pl} corresponds to a familiar operation in the imprecise probability literature, namely, the {\em probability-to-possibility transform} \citep[e.g.,][]{dubois.etal.2004}. 

The possibility contour in \eqref{eq:pl} extends to a full-blown possibility measure for $\Theta$:
\begin{equation}
\label{eq:poss}
\uPi_x(A) = \sup_{\theta \in A} \pi_x(\theta), \quad A \subseteq \TT. 
\end{equation}
There is also a dual {\em necessity measure} defined via conjugacy, i.e., $\lPi_x(A) = 1 - \uPi_x(A^c)$, and it's easy to show that $\lPi_x(A) \leq \uPi_x(A)$ for all $A \subseteq \TT$ and all $x \in \XX$.  Possibilistic reasoning proceeds exactly as described above.  The difference here compared to at the end of Section~\ref{S:likelihood} is that now there's a universal possibility scale, so it's easy for the user to decide what it means for $\uPi_x(A)$ to be ``small''---or, equivalently, what it means for $\lPi_x(A)$ to be ``large''---and to understand the methodological implications of this decision. 

\begin{thm0}
The IM determined by the possibility contour in \eqref{eq:pl} is {\em (strongly) valid}, i.e., 
\begin{equation}
\label{eq:valid}
\sup_{\Theta \in \TT} \prob_\Theta\{ \pi_X(\Theta) \leq \alpha\} \leq \alpha, \quad \text{for all $\alpha \in [0,1]$}, 
\end{equation}
and, consequently, the possibility measure \eqref{eq:poss} satisfies 
\begin{equation}
\label{eq:valid.pm}
\sup_{\Theta \in \TT} \prob_\Theta\{ \uPi_X(A) \leq \alpha \text{ for some $A \ni \Theta$} \} \leq \alpha, \quad \text{for all $\alpha \in [0,1]$}.
\end{equation}
\end{thm0}

\begin{proof}
More general results are covered in \citet{martin.partial2}.  Claim \eqref{eq:valid} can be verified directly via the aforementioned connection to the familiar relative likelihood-based p-values.  Claim \eqref{eq:valid.pm}---that calibration holds uniformly over all true hypotheses---follows from \eqref{eq:valid} and the fact that $\sup_{\theta \in A} \pi_X(\theta) \leq \alpha$ for some $A \ni \Theta$ if and only if $\pi_X(\Theta) \leq \alpha$. 
\end{proof}

One important consequence of the above theorem is that the possibilistic IM is not afflicted by false confidence the way probabilistic inference is \citep{martin.nonadditive, imchar}.  Specifically, false confidence would arise if the IM tended to assign large $\lPi_X$-values to false hypotheses.   But it follows immediately from \eqref{eq:valid.pm} and the definition of $\lPi_x$ that 
\[ \sup_{\Theta \in \TT} \prob_\Theta\{ \lPi_X(A) \geq 1-\alpha \text{ for some $A \not\ni \Theta$} \} \leq \alpha, \quad \alpha \in [0,1]. \]
That is, the event where {\em any} false hypothesis is assigned a relatively large $\lPi_X$-value has relatively small probability---hence no false confidence.  

The following corollary establishes that the same IM output that can be used for reliable in-sample possibilistic reasoning can also be used to construct statistical methods or procedures that ``work'' in the out-of-sample sense described above. 

\begin{cor0}
Hypothesis testing and confidence set procedures derived from the IM defined above control frequentist error rates at the nominal levels.  That is:
\begin{itemize}
\item For a hypothesis $H: \Theta \in A$, the test that rejects $H$ if and only if $\uPi_X(A) \leq \alpha$ has Type~I error probability bounded above by $\alpha$, and
\vspace{-2mm}
\item The set $C_\alpha(X) = \{\theta: \pi_X(\theta) > \alpha\}$ has coverage probability bounded below by $1-\alpha$. 
\end{itemize} 
\end{cor0}

These are the usual frequentist properties expected of hypothesis tests and set estimators and they follow immediately from \eqref{eq:valid} without any conditions on the models involved, the sample size, etc.  The key role played by the likelihood suggests that these IM-driven methods would be efficient, and it is often the case in applications that they agree with the optimal or best-available methods.  Moreover, the particular result in \eqref{eq:valid.pm} implies much stronger error rate control than the first part of the above corollary lets on.  Indeed, the error rate control is actually {\em uniform} in the hypotheses $A$, which has certain probativeness consequences \`a la  \citet{mayo.book.2018}; see \citet{cella.martin.probing}.  

Next I'll show four quick illustrations of the IM formulation, using some examples that were interesting to Basu.  The first pair are problematic ``non-regular'' examples in which the minimal sufficient statistic has dimension greater than that of the unknown parameter, as considered in \citep{basu1964, basu1967} and elsewhere.  

\begin{figure}[t]
\begin{center}
\subfigure[Example~\ref{ex:unif.endpoints}, $(x_{(1)},x_{(n)}) = (281.1, 9689.7)$]{\scalebox{0.6}{\includegraphics{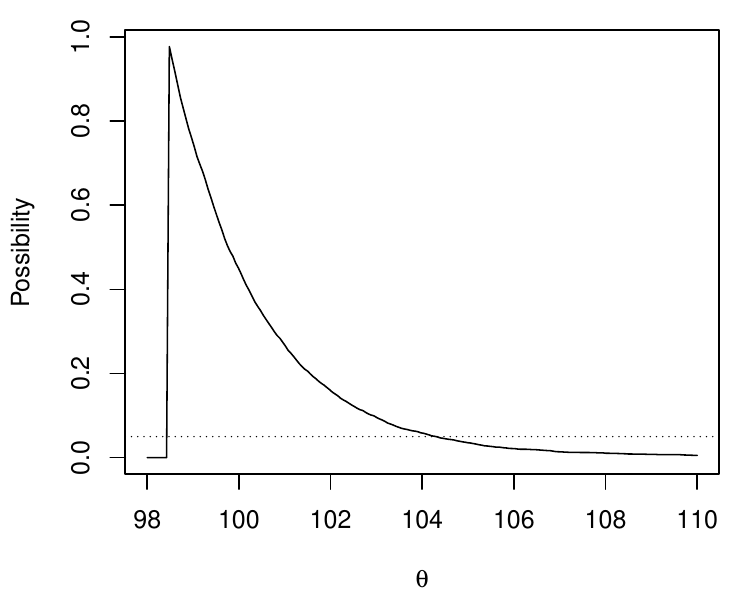}}}
\subfigure[Example~\ref{ex:bvn}, $\hat\theta_x = 0.788$]{\scalebox{0.6}{\includegraphics{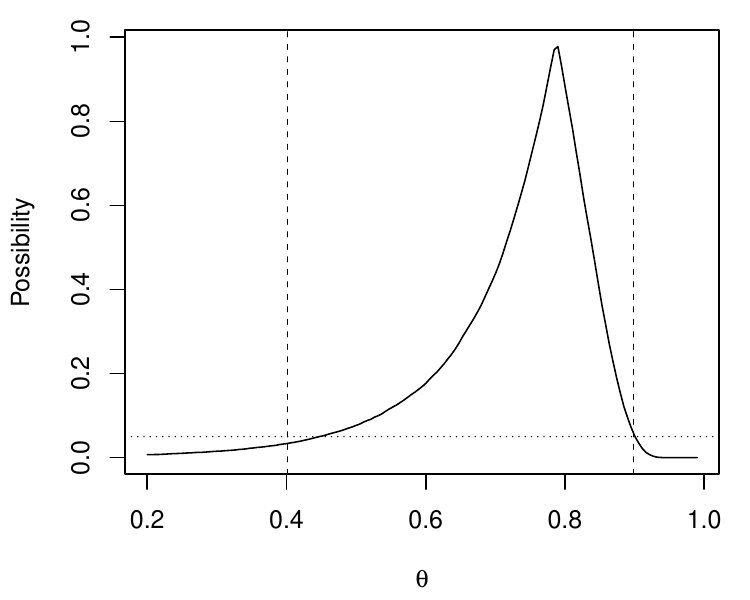}}}
\end{center}
\caption{Plots of the possibility contour functions for Examples~\ref{ex:unif.endpoints}--\ref{ex:bvn}. In Panel~(b), the dashed vertical lines mark the endpoints of an asymptotically efficient 95\% confidence interval based on the ``$r^\star$'' approximation described in \citet{reid2003} and elsewhere.}
\label{fig:uniform1}
\end{figure}

\begin{ex}
\label{ex:unif.endpoints}
Let $X=(X_1,\ldots,X_n)$ consist of iid $\unif\{a(\Theta), a(\Theta) + b(\Theta)\}$ random variables, where $\Theta$ is an unknown scalar but $a(\cdot)$ and $b(\cdot)$ are known functions.  One of Basu's favorites is $a(\theta)=\theta$ and $b(\theta) \equiv 1$, so that $\theta$ is a location parameter.  That's a special group-invariant case, as studied in, e.g., \citet{basu.ghosh.1969}, and the connection between the proposed IM solution and the fiducial/default-prior Bayes solution is presented in \citet{im.group.fiducial}.  Here I consider the model $\prob_\Theta = \unif(\Theta, \Theta^2)$, with unknown $\Theta \in \TT = (1, \infty)$, with endpoint functions $a(\theta) = \theta$ and $b(\theta) = \theta^2 - \theta$.  This problem is ``non-regular'' in the sense that the minimal sufficient statistic---$(X_{(1)}, X_{(n)})$, the extreme order statistics---is two-dimensional while $\Theta$ is a scalar.  Despite this non-regularity, the maximum likelihood estimator is easy to get, $\hat\theta_x = x_{(n)}^{1/2}$, so the relative likelihood is 
\[ \eta(x,\theta) = \frac{L_x(\theta)}{L_x(\hat\theta_x)} = \Bigl\{ \frac{x_{(n)} - x_{(n)}^{1/2}}{\theta^2-\theta} \Bigr\}^n \cdot 1\{x_{(1)} \geq \theta, \, x_{(n)} - x_{(1)} \leq \theta^2-\theta\}, \quad \theta > 1. \]
The corresponding possibility contour based on \eqref{eq:pl} is 
\[ \pi_x(\theta) = \prob_\theta\{ X_{(n)} - X_{(n)}^{1/2} \leq x_{(n)} - x_{(n)}^{1/2}\}, \quad \theta > x_{(1)}. \]
There is no closed-form expression for this, but it's easy to approximate via Monte Carlo: given $\Theta=\theta$, $X_{(n)}$ is distributed as $\theta + (\theta^2-\theta) \bet(n,1)$, a transformed beta random variable.  For an illustration, I follow \citet[][Ex.~4]{hannig.review} and consider data of size $n=25$ with $(x_{(1)}, x_{(n)})=(281.1, 9689.7)$, so that the maximum likelihood estimator is $\hat\theta_x \approx 98.4$. Figure~\ref{fig:uniform1}(a) shows the possibility contour and, from this, we can readily evaluate $\uPi_x(A)$ for any $A$ and/or extract a $100(1-\alpha)$\% confidence interval by thresholding the contour at level $\alpha$.  This IM agrees with that developed in \citet{imunif} and shown there to be both valid and as efficient as existing solutions.  
\end{ex}

\begin{ex}
\label{ex:bvn}
Suppose that, $X=(X_1,\ldots,X_n)$ consists of $n$ iid random variable pairs $X_i = (X_{i1}, X_{i2})$ which are bivariate normal with mean 0, variance 1, and unknown correlation $\Theta \in \TT=[-1,1]$.  It's easy to check that the minimal sufficient statistic is two-dimensional, compared to one-dimensional $\Theta$, so there's non-regularity like in Example~\ref{ex:unif.endpoints}.  Consequently, inference based on the sampling distribution of, say, the maximum likelihood estimator would be inefficient due to a loss of information.  As \citet{basu1964} showed, there's no guidance on what ancillary statistic one should condition on to recover the lost information---both partial data sets $(X_{i1},\ldots,X_{n1})$ and $(X_{i2},\ldots,X_{n2})$ are equally good ancillary statistics---so it's not clear how to proceed.  For this reason, care is needed in developing valid and efficient methods for inference on $\Theta$.  

There's no closed-form expression for the maximum likelihood estimator and, consequently, there's no closed-form expression for the relative likelihood or the proposed IM's possibility contour $\theta \mapsto \pi_x(\theta)$ in \eqref{eq:pl}.  Fortunately, it's not too difficult to perform the required computations numerically (e.g., root-finding and Monte Carlo), and I have done so for a simulated data set of size $n=10$ with $\Theta=0.75$.  Figure~\ref{fig:uniform1}(b) shows the possibility contour plot; for reference, the maximum likelihood estimator in this case is $\hat\theta_x = 0.788$.  One can easily read off a 95\% confidence interval by thresholding the contour at level $\alpha=0.05$.  For comparison, the vertical lines mark the endpoints of the 95\% confidence interval determined by the efficient ``$r^\star$'' asymptotic approximation \citep[e.g.,][]{brazzale.davison.2008, reid2003}.
The difference between these and the IM's confidence limits is negligible, so the IM solution is exactly, provably valid and efficient. 
\end{ex}

The second pair of illustrations are closely tied to Basu's interests in and deep insights concerning finite-population sampling \citep[e.g.,][]{basu1969, basu1971}. 

\begin{figure}[t]
\begin{center}
\subfigure[Example~\ref{ex:unif.discrete}, $x_{(n)}=5$]{\scalebox{0.6}{\includegraphics{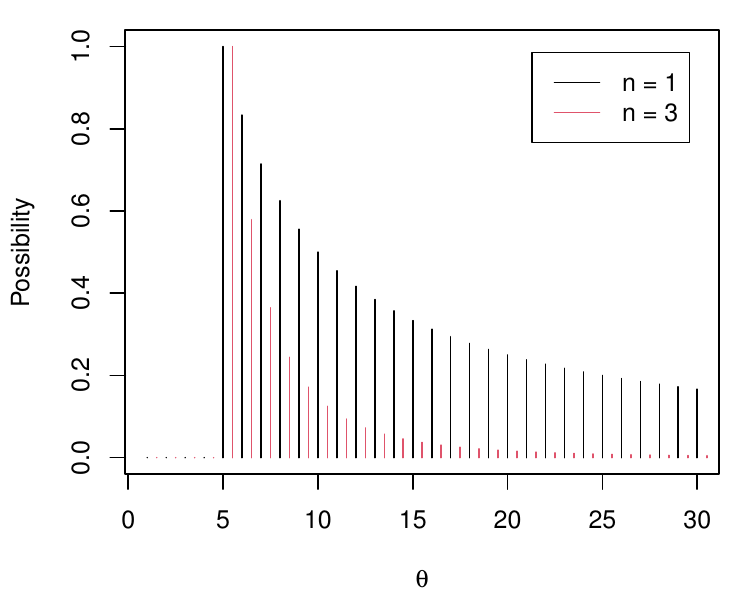}}}
\subfigure[Example~\ref{ex:urn}, $x=1$]{\scalebox{0.6}{\includegraphics{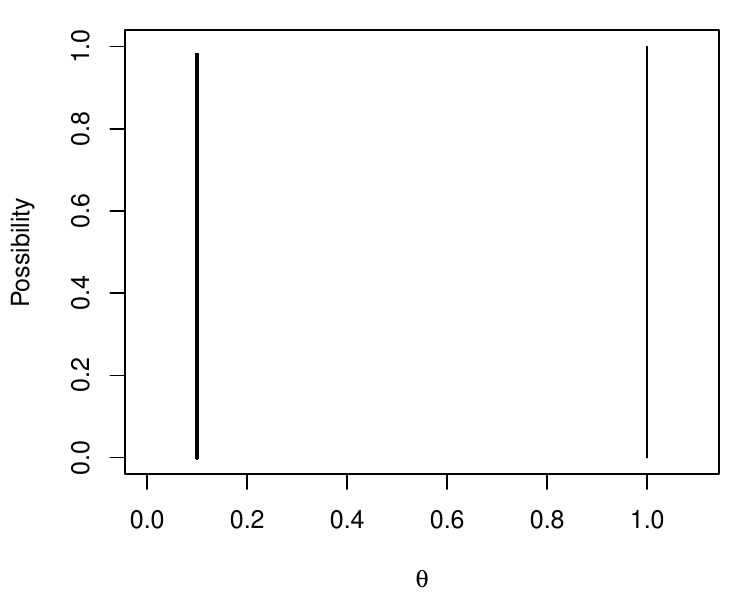}}}
\end{center}
\caption{Possibility contour plots for Examples~\ref{ex:unif.discrete}--\ref{ex:urn}. In Panel~(a), the red spikes have been shifted to the right so they don't overlap with the black. In Panel~(b), the thick vertical line is made up of 980 spikes of height 0.98 around the value $\theta=0.1$.}
\label{fig:uniform2}
\end{figure}

\begin{ex}
\label{ex:unif.discrete}
Suppose $X=(X_1,\ldots,X_n)$ is an iid sample from $\unif\{1,2,\ldots,\Theta\}$, where $\Theta$ is an unknown natural number.  In this case, the likelihood function is 
\[ L_x(\theta) = \theta^{-n} \cdot 1(\theta \geq x_{(n)}), \quad  \theta=1,2,\ldots \]
and the maximum likelihood estimator is $\hat\theta_x = x_{(n)}$, so it's easy to show that the IM's possibility contour based on \eqref{eq:pl} is  
\[ \pi_x(\theta) = ( x_{(n)} /\theta )^n, \quad \theta=x_{(n)}, x_{(n)} + 1, \ldots \]
A plot of this contour function is shown in Figure~\ref{fig:uniform2}(a) based on $x_{(n)}=5$ with two values of $n$; the vertical spikes emphasize that it's a function defined only on the integer values.  Clearly, these are not probability masses since they don't sum to 1.  The maximum possibility value of 1 is attained at $\theta=x_{(n)}$, decreasing thereafter, and the extended possibility measure \eqref{eq:poss} on general hypotheses can be readily evaluated as needed.  Note that the possibility contour vanishes much more rapidly for $n=3$ compared to $n=1$, which is sign of the efficiency gain with a larger sample size. 
\end{ex}

\begin{ex}
\label{ex:urn}
Consider the example in \citet[][p.~240]{basu1975} involving an urn that contains 1000 balls: 20 are labeled with $\Theta$ and the remaining 980 are labeled with the values $a_1 \Theta,\ldots,a_{980} \Theta$, where the $a_j$'s are distinct known values in the interval $[9.9, 10.1]$.  Let $X$ denote the value on a single randomly chosen ball from this urn.  The likelihood is 
\[ L_x(\theta) = \begin{cases} 0.02 & \text{if $\theta=x$} \\ 0.001 & \text{if $\theta \in \{a_1^{-1} x, \ldots, a_{980}^{-1} x\}$} \\ 0 & \text{otherwise}. \end{cases} \]
Basu designed this example to highlight some unusual behavior of the maximum likelihood estimator, in particular, that $\hat\theta_X=X$ is far from $\Theta$ with $\prob_\Theta$-probability 0.98.  From here it's not difficult to show that the IM's possibility contour is 
\[ \pi_x(\theta) = \begin{cases} 1 & \text{if $\theta=x$} \\ 0.98 & \text{if $\theta \in \{a_1^{-1} x, \ldots, a_{980}^{-1} x\}$} \\ 0 & \text{otherwise}. \end{cases} \]
A plot of this is shown in Figure~\ref{fig:uniform2}(b) with $x=1$.  The aforementioned unusual behavior of the maximum likelihood estimator is not an issue here because there's no compelling reason to single out $\hat\theta_x$ when all the other values---which are very close to $\Theta$ when $\hat\theta_x$ isn't---are similarly highly plausible.  
\end{ex}


To conclude, the framework that I'm proposing here is a viable candidate for Basu's {\em via media}.  It combines the (Bayesian-like) in-sample possibilistic reasoning with the (frequentist-like) calibration that guarantees the derived methods ``work'' in the out-of-sample sense that's relevant to users of statistical methods.  This combination can't be achieved without venturing into imprecise probability territory.

\section{Valid IMs and the likelihood principle}
\label{S:favorites}




The IM framework I put forward in Section~\ref{S:new} doesn't satisfy the likelihood principle.  That is, despite being largely relative likelihood-driven, the possibility contour isn't {\em fully} determined by the relative likelihood---it depends on the model $\{\prob_\theta: \theta \in \TT\}$, on the sample space $\XX$, etc.~via the probability calculation in \eqref{eq:pl}---so it fails to satisfy the likelihood principle.  But this doesn't mean that it's impossible to achieve the likelihood principle (or something close enough to it), if desired, through some adjustments.  Remember, we're after a {\em via media}, so certain trade-offs should be expected to meet today's methods-focused needs.  These adjustments will also highlight the flexibility that an imprecise-probabilistic framework affords the statistician. 

Recall that the posited model $\{\prob_\theta: \theta \in \TT\}$ and observed data $X=x$ determine a relative likelihood $\eta(x,\cdot)$, but not uniquely.  That is, in general there's a class of models that all produce the same $\eta(x,\cdot)$ for (almost) all $x$.  To give the reader some context, what I have in mind are different data-collection procedures, e.g., sampling designs, stopping rules, etc., for investigating the same scientific question.  To ensure that these developments make sense mathematically, I'll reinterpret the data $X$ as whatever's needed to determine that relative likelihood.  Given $\eta: \XX \times \TT \to [0,1]$, let 
\[ \rlik^\star = \rlik^\star(\eta,\TT) = \{\prob_\theta^{(m)}: \theta \in \TT, \, m \in \MM^\star\}, \]
where $m \in \MM^\star$ is a generic model index, denote the collection of {\em all} probability distributions on $\XX$, parametrized by $\theta \in \TT$, with density/mass function $p_\theta^{(m)}(x)$ that satisfies 
\[ \frac{p_\theta^{(m)}(x)}{\sup_\vartheta p_\vartheta^{(m)}(x)} = \eta(x,\theta), \quad \text{for all $\theta \in \TT$, $m \in \MM^\star$, and (almost) all $x \in \XX$}. \]
That is, any one of these candidate models determines $\eta$, and then the data analyst collects in $\MM^\star$ all the models with equivalent relative likelihood. 

For a concrete example, consider a sequence of Bernoulli trials where the data is a pair consisting of the number of trials performed and the number of successes observed; write this as $x=(n,y)$, where $n$ is the number of trials and $y$ is the number of successes.  There are, of course, a variety of models for data of this type, depending on how the experiment is performed.  If $n$ is fixed in advance, then $y$ would be considered ``data,'' and a binomial model would be appropriate.  Alternatively, if  $y$ is fixed in advance, then $n$ is the ``data,'' and a negative binomial model would be appropriate.  As is well-known, both of these have relative likelihood 
\[ \eta(x, \theta) = \Bigl( \frac{n \theta}{y} \Bigr)^y \Bigl( \frac{n-n\theta}{n-y} \Bigr)^{n-y}, \quad \theta \in [0,1], \quad x=(n,y). \]
While the above two designs might be the most common in practice, these aren't the only two models in $\MM^\star$ for $\eta$ as above; there are many more, one for each proper stopping rule.  Example~21 of \citet{bergerwolpert1984} offers a setup wherein $x=(n,y)$ can take one of three possible values, namely, $(1,1)$, $(2,0)$, or $(2,1)$, i.e., stop the study after the first trial if it's a success, otherwise stop after the second trial.  

The data analyst might be able to eliminate some of the equivalent models so, in general, consider a sub-collection $\MM \subseteq \MM^\star$.  In the Bernoulli trial illustration above, the data analyst might not know what stopping rule was used, but if he knows that some {\em weren't} used, then those can be omitted from $\MM$.  The embellishment I'm suggesting here, natural from an imprecise-probabilistic point of view, is to define a new possibility contour by maximizing the right-hand side of \eqref{eq:pl} over models: 
\begin{equation}
\label{eq:pl.lp}
\pi_x(\theta \mid \MM) = \sup_{m \in \MM} \prob_\theta^{(m)}\{ \eta(X,\theta) \leq \eta(x,\theta) \}, \quad \theta \in \TT. 
\end{equation}
Since there won't be any chance for confusion in what follows, I'll drop the dependence on $\MM$ in the notation above, and just write ``$\pi_x(\theta)$'' for the right-hand side in \eqref{eq:pl.lp}.  Note that each $\theta \mapsto \prob_\theta^{(m)}\{ \eta(X,\theta) \leq \eta(x,\theta)\}$ takes value 1 at the maximum likelihood estimator, so the right-hand side satisfies $\sup_\theta \pi_x(\theta)=1$, hence is a possibility contour.  Therefore, I can define a possibility measure $\uPi_x(A) = \sup_{\theta \in A} \pi_x(\theta)$ exactly as before, and the same in-sample possibilistic reasoning can be applied.  It's also immediately clear that the IM validity property \eqref{eq:valid} holds here too, so the derived methods are provably reliable.  But the validity conclusions are broader because they hold uniformly over the models in $\rlik$.  The broader validity conclusions come at a price though: the supremum over $\MM$ implies that the possibility contour in \eqref{eq:pl.lp} is no more tightly concentrated than that corresponding to any $m$-specific model, hence a potential loss of efficiency, e.g., larger confidence sets.  This loss of efficiency is unavoidable if one wants the likelihood principle and reliability guarantees.  In any case, what I'm proposing is very much {\em via media} in spirit since the practitioner can control how close he is to satisfying the likelihood principle---and how much efficiency he stands to lose---through his choice of $\MM \subseteq \MM^\star$.  

\begin{figure}[t]
\begin{center}
\subfigure[$x=(n,y)=(10, 3)$]{\scalebox{0.6}{\includegraphics{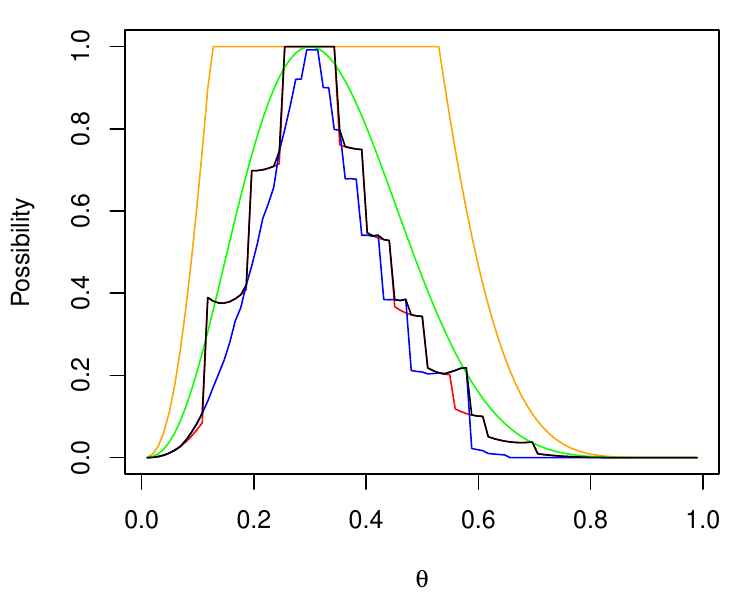}}}
\subfigure[$x=(n,y)=(16, 11)$]{\scalebox{0.6}{\includegraphics{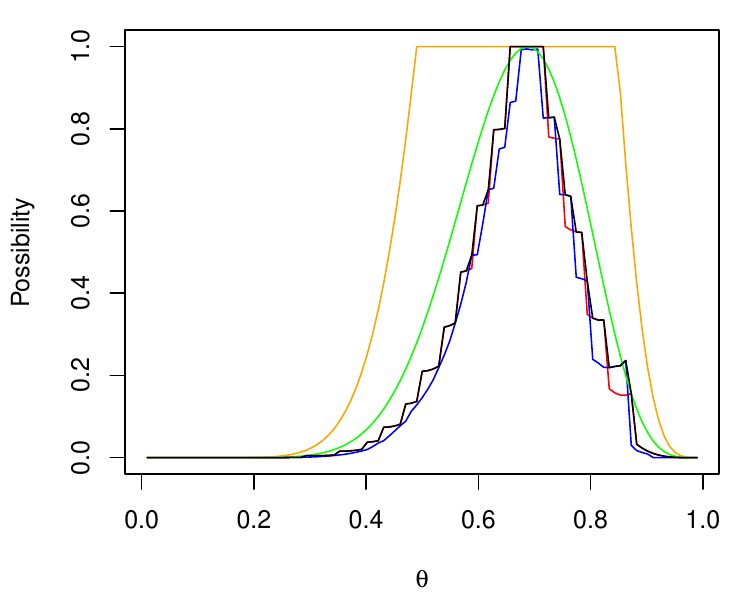}}}
\end{center}
\caption{Plots of the model-specific possibility contours (binomial is red, negative binomial is blue) and the pointwise maximum (black) in \eqref{eq:pl.lp}. The green curve is the relative likelihood $\eta$ and the orange curve is the truncated $\eta_Q$ in \eqref{eq:e.poss}, with $Q = \unif(0,1)$.}
\label{fig:bernoulli}
\end{figure}

Returning to the Bernoulli trial illustration, suppose that $\MM$ contains just the binomial and negative binomial models.  Figure~\ref{fig:bernoulli} shows plots of the two model-specific possibility contours as in \eqref{eq:pl} and the combined version in \eqref{eq:pl.lp} for two different data sets $x=(n,y)$.  This plot highlights the point that, thanks to sharing the same $\eta$, the model-specific possibility contours have overall similar shapes.  This means that the pointwise maximum in \eqref{eq:pl.lp} isn't going to be too much different from the individual curves, which is apparent in the plots.  So, for example, the confidence intervals obtained by thresholding the three curves at level $\alpha$ are all about the same.  The difference is in which model(s) the coverage probability claims apply to: the interval determined by \eqref{eq:pl.lp} satisfies the coverage probability claim for {\em both} the binomial and negative binomial models.   

The general case in \eqref{eq:pl.lp} is computationally intimidating, and I don't presently have any recommendations on how this can be carried out efficiently, but bounds may be available; see below.  I imagine, however, that a data analyst who is seriously concerned about both reliability and satisfying the likelihood principle can identify a relatively small finite set of plausible models in $\MM^\star$ that deserve consideration. Then the computations wouldn't be much more difficult than those needed to generate the plots in Figure~\ref{fig:bernoulli}.  It's in the user's best interest, after all, to be parsimonious in their choice of $\MM$, since an overly generous choice will lead to unnecessary loss of efficiency.  

I'll conclude this section by discussing the case where $\MM = \MM^\star$, i.e., where the user is entertaining literally all the models that share a common relative likelihood.  \citet{walley2002}, for instance, develops a framework of (imprecise) probabilistic inference that both satisfies the likelihood principle and achieves a version of the validity result in \eqref{eq:valid}.  Let $Q$ denote a generic prior probability distribution for $\Theta$ on $\TT$ and define 
\[ \eta_Q(x,\theta) = \frac{L_x(\theta)}{\int L_x(\vartheta) \, Q(d\vartheta)}, \quad \theta \in \TT. \]
If $\theta$ were some specific parameter value $\theta_0$, then $\eta_Q(x, \theta_0)$ might be referred to as the {\em Bayes factor} for testing the null hypothesis ``$\Theta = \theta_0$'' against the alternative hypothesis ``$\Theta \sim Q$.''  In addition to the Bayesian interpretation, $\eta_Q$ has a few interesting and relevant properties.  First, since $\eta_Q$ only depends on the likelihood function up to proportionality, all the models in $\MM^\star$ yield the same $\eta_Q$ for a given $Q$, just like with $\eta$.  Second, the reciprocal of $\eta_Q$ (but not of $\eta$) satisfies 
\begin{equation}
\label{eq:evalue}
\E_\theta^{(m)} \{\eta_Q(X,\theta)^{-1}\} = 1, \quad \text{for all $\theta \in \TT$, $m \in \MM$}, 
\end{equation}
where $\E_\theta^{(m)}$ denotes expected value with respect to $\prob_\theta^{(m)}$.  This follows easily because $x \mapsto \int p_\vartheta(x) \, Q(d\vartheta)$ defines a density/mass function.  This shows that $\eta_Q^{-1}$ determines an {\em e-value} or {\em e-process} \citep[e.g.,][]{evalues.review}.  An immediate consequence is that $\eta_Q$ achieves a property similar like in \eqref{eq:valid}: by Markov's inequality and \eqref{eq:evalue}
\begin{equation}
\label{eq:markov}
\prob_\theta^{(m)}\{ \eta_Q(X,\theta) \leq \alpha\} \leq \alpha \, \E_\theta^{(m)}\{ \eta_Q(X,\theta)^{-1} \} = \alpha, \quad \text{$\alpha \in [0,1]$, $\theta \in \TT$, $m \in \MM$}. 
\end{equation}
Therefore, $\eta_Q$ can be readily used to construct valid tests and confidence intervals as in the above corollary.  Since the ``$\leq \alpha$'' bound in \eqref{eq:markov} holds uniformly in $m \in \MM$, the conclusions can be strengthened to ``anytime valid'' \citep{evalues.review}, but I'll not explain these details here.  Note, however, that $\theta \mapsto \eta_Q(x,\theta)$ is not a possibility contour; since $\eta_Q(X,\theta)^{-1}$ has expected value 1, aside from trivial cases, it'll surely take values greater than 1.  But it can be truncated to a possibility contour, 
\begin{equation}
\label{eq:e.poss}
\theta \mapsto \eta_Q(x,\theta) \wedge 1, \quad \theta \in \TT. 
\end{equation}
One can imagine, however, that the procedures derived from thresholding the $Q$-specific possibility contour in \eqref{eq:e.poss} would be conservative, since the validity guarantees would have to hold for any user-specified $Q$.  This conservatism is apparent in Figure~\ref{fig:bernoulli}.  

The relative likelihood $\eta$ and the corresponding $\eta_Q$ are related via 
\[ \eta(x,\theta) = \inf_Q \eta_Q(x,\theta), \]
where the infimum is over all probability measures on $\TT$, and it's attained at a measure that assigns probability~1 to set of maximizers of the likelihood for the given data $x$.  Then the following strategy is tempting: first define a $Q$-specific possibility contour 
\[ \pi_x(\theta \mid Q) = \pi_x(\theta \mid \MM, Q) = \sup_{m \in \MM} \prob_\theta^{(m)}\{ \eta_Q(X,\theta) \leq \eta_Q(x,\theta)\}, \quad \theta \in \TT, \]
and then try removing the dependence on $Q$ by optimizing again, i.e., 
\[ \tilde\pi_x(\theta) = \inf_Q \pi_x(\theta \mid Q), \quad \theta \in \TT. \]
This satisfies the likelihood principle, since it doesn't depend on any particular model $m$ in $\MM=\MM^\star$.  Moreover, by the bound in \eqref{eq:markov}, 
\[ \tilde\pi_x(\theta) \leq \inf_Q \{ \eta_Q(x,\theta) \wedge 1\} = \eta(x,\theta), \quad \theta \in \TT. \]
Remember, the relative likelihood on the right-hand side above is a possibility contour but the corresponding possibilistic IM isn't valid---e.g., it doesn't satisfy \eqref{eq:evalue} because $x \mapsto \sup_\vartheta p_\vartheta(x)$ isn't a density/mass function.  The problem is that $\eta$ tends to be too small which, together with the above bound, implies that $\tilde\pi_x$ doesn't define a valid IM either.  We do, however, get the following insights:
\begin{itemize}
\item the relative likelihood-based possibility contour in \eqref{eq:pl.lp}, with $\MM = \MM^\star$, will tend to be less tightly concentrated than the relative likelihood itself, and 
\vspace{-2mm}
\item at least intuitively, the relative likelihood-based possibility contour in \eqref{eq:pl.lp} ought to be more tightly concentrated then $\eta_Q \wedge 1$ for any particular $Q$. 
\end{itemize} 
These observations are apparent in Figure~\ref{fig:bernoulli}.  An interesting open question is whether the vague notion of ``tight'' that I'm using above could be related to the familiar, well-defined notion of {\em specificity} in the possibility theory literature \citep[e.g.,][]{dubois.prade.1986}.  In any case, I'd feel comfortable upper bounding the possibility contour \eqref{eq:pl.lp} in the challenging case with $\MM = \MM^\star$ by $\eta_Q \wedge 1$ for some not-too-tightly-concentrated $Q$.



\section{Conclusion}
\label{S:discuss}

In this paper, I revisited the {\em possibility} of achieving a middle-ground between Basu's Bayesian and frequentist poles.  Resolving this open question would go a long way towards pinpointing statisticians' contribution and securing our seat at the data science table.  While Fisher's efforts fell short, my claim is that there's still hope.  The key new observation, as I described in Section~\ref{S:likelihood}, is that likelihood (model + data) is insufficient to reliably support probabilistic inference.  This helps to justify consideration of other non-traditional modes of uncertainty quantification.  Furthermore, I've argued (here and elsewhere) that likelihod can reliably support possibilistic inference, and I've offered a framework in which this can be carried out.  There's still some work to be done, but I think almost all of the relevant details have been worked out in \citet{martin.partial2}.  If I'm wrong and this isn't the {\em via media} that Fisher and others have been looking for, then I urge the reader to reach out and let me know what I'm missing.  

For further developments, I'm very excited about the potential for incorporating partial prior information into the possibilistic IM, like I mentioned briefly in Section~\ref{S:new}.  A practically important and challenging problem---another favorite of Basu's \citep[e.g.,][]{basu1977, basu1978}---that tends to get overlooked is marginal inference in the presence of nuisance parameters.  The possibility-theoretic framework offers a straightforward marginalization procedure that preserves validity; this is via the {\em extension principle} of \citet{zadeh1975a}.  The downside is that this straightforward marginalization tends to be inefficient.  To avoid this loss of efficiency, some form of dimension reduction is needed.  The general profiling strategy I proposed in \citet[][Sec.~7]{martin.partial2} seems promising, but I've since realized that, in certain cases, more efficient marginal inference can be achieved using other strategies besides profiling.  So there are still more insights to be gleaned from Basu on this important question, and I'll report on these details elsewhere.

\section*{Acknowledgments}

The author thanks an anonymous reviewer for valuable feedback on a previous version of the manuscript.  This work is partially supported by the U.S.~National Science Foundation, grant SES--2051225.

%
%

\bibliographystyle{apalike}
\bibliography{/Users/rgmarti3/Dropbox/Research/mybib}

\end{document}